\def\mycmd{1} 
\if\mycmd0
\documentclass{article} 
\else
\documentclass{elsarticle} 
\fi

\usepackage{graphicx,subcaption,pgfplots} 
\usepackage[utf8]{inputenc}
\usepackage{todonotes,booktabs}
\usepackage{fullpage}
\usepackage{amsmath, amssymb}
\usepackage{bm}
\usepackage{isomath,pifont}
\usepackage{amsthm} 
\usepackage{mathtools} 
\if\mycmd0
\usepackage{authblk}
\usepackage{mathpazo}
\fi
\renewcommand*{\vec}{\vectorsym}
\newcommand{\ten}{\tensorsym}

\newcommand{\cmark}{\ding{51}}%
\newcommand{\xmark}{\ding{55}}%

\DeclareMathOperator{\dive}{\text{div}}
\DeclareMathOperator{\grad}{\nabla}
\newcommand{\der}[2]{\frac{d#1}{d#2}}
\newcommand{\parder}[2]{\frac{\partial #1}{\partial #2}}

\newcommand{\R}{\mathbb{R}}

\newcommand{\porosityL}{\Phi} 
\newcommand{\porosityE}{\phi} 

\newtheorem{remark}{Remark}

\newcommand*{\TitleFont}{%
      \usefont{\encodingdefault}{\rmdefault}{b}{n}%
      \fontsize{16}{20}%
      \selectfont}

\begin{document}
\if\mycmd0
\title{\TitleFont \textbf{Efficient reference configuration formulation in fully nonlinear poroelastic media}}

\author[1,*]{Nicolás A. Barnafi}

\author[2]{Argyrios Petras}

\author[2,3]{Luca Gerardo-Giorda}

\affil[1]{Center for Mathematical Modelling (CMM), Santiago, Chile} 
\affil[2]{Johann Radon Institute for Computational and Applied Mathematics (RICAM), Austrian Academy of Sciences, Linz, Austria} 
\affil[3]{Institute for Mathematical Methods in Medicine and Data-Based Modelling, Johannes Kepler University, Linz, Austria}


\date{}

\maketitle
\begin{abstract}
    Typical pipelines for model geometry generation in computational biomedicine stem from images, which are usually considered to be at rest, despite the object being in mechanical equilibrium under several forces. We refer to the stress-free geometry computation as the reference configuration problem, and in this work we extend such a formulation to the theory of fully nonlinear poroelastic media. The main steps are (i) writing the equations in terms of the reference porosity and (ii) defining a time dependent problem whose steady state solution is the reference porosity. This problem can be computationally challenging as it can require several hundreds of iterations to converge, so we propose the use of Anderson acceleration to speed up this procedure. Our evidence shows that this strategy can reduce the number of iterations up to 80\%. In addition, we note that a primal formulation of the nonlinear mass conservation equations is not consistent due to the presence of second order derivatives of the displacement, which we alleviate through adequate mixed formulations. All claims are validated through numerical simulations in both idealized and realistic scenarios.
\end{abstract}
\else
\begin{frontmatter}

 \title{Efficient reference configuration formulation in fully nonlinear poroelastic media}

 \author[1,*]{Nicolás A. Barnafi}

 \author[2]{Argyrios Petras}

 \author[2,3]{Luca Gerardo-Giorda}

\address[1]{Center for Mathematical Modelling (CMM), Universidad de Chile, Beauchef 851, Santiago, Chile. \texttt{nbarnafi@cmm.uchile.cl}} 
\address[2]{Johann Radon Institute for Computational and Applied Mathematics (RICAM), Austrian Academy of Sciences, Linz, Austria} 
\address[3]{Institute for Mathematical Methods in Medicine and Data-Based Modelling, Johannes Kepler University, Linz, Austria}

\date{\today}
 
\begin{abstract}
Typical pipelines for model geometry generation in computational biomedicine stem from images, which are usually considered to be at rest, despite the object being in mechanical equilibrium under several forces. We refer to the stress-free geometry computation as the reference configuration problem, and in this work we extend such a formulation to the theory of fully nonlinear poroelastic media. The main steps are (i) writing the equations in terms of the reference porosity and (ii) defining a time dependent problem whose steady state solution is the reference porosity. This problem can be computationally challenging as it can require several hundreds of iterations to converge, so we propose the use of Anderson acceleration to speed up this procedure. Our evidence shows that this strategy can reduce the number of iterations up to 80\%. In addition, we note that a primal formulation of the nonlinear mass conservation equations is not consistent due to the presence of second order derivatives of the displacement, which we alleviate through adequate mixed formulations. All claims are validated through numerical simulations in both idealized and realistic scenarios.
\end{abstract}

\begin{keyword}
Nonlinear poroelasticity \sep Anderson acceleration \sep Reference configuration \sep Inverse design  \sep Robust discretizations \sep Mixed formulations



\end{keyword}

\end{frontmatter}
\fi

\section{Introduction}
Porous media, understood as objects that exhibit a solid phase and an interconnected network through which a fluid flows, are ubiquitous in nature and can be found in materials like wood, rubber, and stones, in animals such as sponges, and in organs like lungs and hearts. See \cite{bear1987theory} for further examples and applications. These objects may deform, but they preserve their porous structure under the action of forces. These forces can be either external or internal from within the pores, and in fact internal pressures in organs play a significant role in their function.  For instance, when a baby takes their first breath, it generates a pressure gradient that pushes the fluid into the interstitial tissue, allowing it to be removed by pulmonary and lymphatic circulations \cite{LoMauro65}. 
The porous structure of biological tissues is essential not only for the normal function of organs through blood perfusion, but also for drug delivery and pharmacological treatment of various diseases \cite{algranati2010mechanisms}. 

Biomedical models typically consider biological tissues as elastic continua, disregarding their complex vascular network that provides blood, which is essential for their function. However, the extreme complexity of vascular networks has motivated the use of a porous media modeling approach in which only the average macroscopic behavior of flow in such networks is described \cite{huyghe1992porous}. Most models use multiple compartments to account for the blood vessels (arteries, veins, and capillaries), the solid part, and even the interstitial space of tissues undergoing large deformations, which applies mainly to the lungs \cite{patte2022quasi} and the heart \cite{cookson2012novel, chapelle2010poroelastic,barnafi2022perfusion,barnafi2022edema}. 

Biological tissues and organs are subject to external forces (gravity, rib cage) and internal forces, which can be either given by air pressure (lungs) or fluid pressure (vessels, lymphatic system, extracellular space), among others. These forces are already present when a patient's geometry is computed from medical imaging, so it becomes fundamental to find a reference configuration whose deformed state is the one computed. This problem, known as the inverse design problem \cite{govindjee1996computational} or prestress problem \cite{gee2009prestressing}, has been successfully applied to the case of a beating (purely elastic) heart \cite{rausch2017augmented, peirlinck2018modular,petras2023mechanoelectric,barnafi2023reconstructing}.

To the best of our knowledge, porous media models using a prestress configuration have only been developed in \cite{patte2022quasi} under strong assumptions that remove the mass conservation law from the system. This work has been motivated by this missing reference configuration model that accommodates to arbitrary nonlinear porous media undergoing large deformations. For this, we propose a novel poroelastic media equation that yields the reference porosity, together with the reference configuration, as an unknown initial condition. 

The resulting model presents several difficulties beyond being extremely nonlinear. The main two are: (i) it can take many time iterations to yield a reference configuration, and (ii) it presents second order derivatives of the displacement--an $\vec H^1$ field--which we refer to as the \emph{primal inconsistency} error. To address these limitations, we propose for (i) to use a fixed-point acceleration technique to reduce the number of iterations, and for (ii) to formulate mixed models that do not suffer from primal inconsistency. The resulting models are accurate and present a framework that is also compatible with their solution in extreme scale computing. We validate our claims through several numerical tests.

This work is structured as follows: in Section~\ref{section:model} we present both the forward equation and the proposed novel reference configuration formulation. In Section~\ref{section:numerical} we present all details pertaining the discretization of the obtained model in order to have a robust and efficient solution strategy. In this section we also present how the computation of a stationary state can be interpreted as a fixed-point problem, and how this can be accelerated. In Section~\ref{section:mixed}, we present our mixed models that avoid the primal consistency error, and comment on their theoretical properties. We present several numerical tests in Section~\ref{section:tests}, and conclude with a discussion of our results in Section~\ref{section:discussion}. 

To avoid confusing notation, throughout this work we will refer to the direct problem as the \emph{forward} problem. Instead, we will refer to the proposed system as the reference configuration model. Also, we will denote scalar, vector, and tensor quantities with regular, bold, and slanted fonts such that $L^2$, $\vec L^2$, and $\ten L^2$ are the corresponding square integrable spaces.

\section{Model formulation}\label{section:model}
In this section we show the formulation of the forward and reference configuration models of nonlinear poroelasticity. We provide all the balance laws related to nonlinear poroelasticity in both Lagrangian and Eulerian coordinates, as both configurations are relevant for devising the forward and reference configuration problems respectively. For further explanation on the model, see \cite{Coussy2004}. 

Consider a reference domain $\Omega_0$ and a displacement field $\vec d:\Omega_0\to \R^3$ such that a reference (Lagrangian) point $\vec X$ is warped into a spatial (current, Eulerian) one $\vec x$ through $\vec x(\vec X) = \vec X + \vec d(\vec X)$. We denote the deformed configuration as $\Omega = \vec x(\Omega_0)$, and further define the strain tensor and its determinant as $\ten F \coloneqq \ten I + \grad_{\vec X}\vec d$ and $J \coloneqq \det \ten F$ respectively. As done with the definition of $\ten F$, whenever relevant we will clarify the coordinates with respect to which we are deriving as $\grad_{\vec X}$ and $\grad_{\vec x}$. The reference and current outwards normal vector will be denoted with $\vec N$ and $\vec n$ respectively.

The geometry $\Omega_0$ is assumed to be a porous medium, where fluid flows through a complex interconnected network structure. Being a continuum theory, both solid and fluid phases coexist spatially. The percentage of the fluid phase in $\Omega$ is given by the Eulerian porosity $\porosityE$, so that naturally the solid porosity is given by $\porosityE_s = 1 - \porosityE$. The Lagrangian porosity $\porosityL$ is defined as $\porosityL = J\porosityE$, such that
        $$ \int_{\omega_0} \porosityL\,dX = \int_\omega \porosityE \,dx, $$
for all $\omega_0\subset\Omega_0$, $\omega=\vec x(\omega_0)$.  Note that all developments and conclusions given in this work can be easily extended to the case of more than one fluid phase, known as \emph{multi-compartment} of \emph{multi-phase} fluids.

\subsection{Forward problem}

In this section, we recapitulate the balance laws for the forward model, given by a quasi-static poroelastic solid. These  are: (i) conservation of linear momentum, (ii) mass conservation, and (iii) solid incompressibility (see \cite{Coussy2004} for their derivations). We highlight that incompressibility is not strictly necessary for this formulation, but we nevertheless decided to include it. This is motivated by our target applications and also by the fact that it can be easily approximated using quasi-incompressibility or even removed, thus rendering the methodology more general.

\paragraph{Equilibrium equation.} Consider a given Helmholtz potential $\Psi$ such that the Piola-Kirchhoff stress tensor is given by $\ten P \coloneqq \parder{\Psi}{\ten F}$. Then, the linear momentum conservation is given by 
    \begin{equation}\label{eq:momentum}
    -\dive \ten P(\ten F) = \vec g \qquad \text{in $\Omega_0$},
    \end{equation}
for any given bulk force $\vec g$, and it represents the balance of forces. Despite the popularity of this formulation, typically forces are known in the deformed configuration (denoted by $\vec g_e$) and not in the reference configuration. Note that performing the pull-back of the force $\vec g_e$ to the reference configuration yields the relationship $\vec g = J\vec g_e$, which stems from the change of variables formula $dx = J dX$. Boundary conditions can be of Dirichlet $\vec d = \bar{\vec{d}}$, Neumann $\ten P\vec N = \bar{\vec t}$, and Robin $\ten P\vec N = -\alpha \vec d$ type in non-overlapping portions of the boundary.

\paragraph{Mass conservation.} For the same potential $\Psi$ considered before, we define the pressure $p$ to be given by $p \coloneqq \parder{\Psi}{\porosityL}$. This fact is deduced from the thermodynamics of poroelastic media \cite{Coussy2004}. The mass conservation of the fluid phase is given by
    \begin{equation}\label{eq:mass}
    \der{\porosityL}{t} +\dive \vec U = 1/{\rho_{f}}\Theta \qquad\text{in $\Omega_0$},
    \end{equation}
where $\Theta$ is a given source/sink term, and $\vec U$ is the relative velocity of the fluid with respect to the solid phase. If the source/sink term is given in the Eulerian configuration, then $\Theta = J\theta$ for a given $\theta$. We adopt a typical choice used in porous media to model the relative fluid velocity $\vec U$, using Darcy's law. This law states that the fluid follows the negative pressure gradient, and is given by
    $$ \vec U \coloneqq -\ten K \grad p,$$
where the tensor $\ten K$ is obtained as the pullback of a given Eulerian permeability tensor $\ten k$, i.e. $\ten K = J\ten F^{-1} \ten k \ten F^{-T}$.  Boundary conditions can be of Dirichlet $\porosityL = \bar{\porosityL}$, Neumann $\ten K\grad p\cdot \vec N = 0$, and Robin type $\ten K\grad p \cdot \vec N = -\alpha \porosityL$. Typically, the porous media is considered to be isolated, meaning that homogeneous Neumann conditions are considered throughout the boundary. The initial condition is given by a known Lagrangian porosity such that 
    $$ \porosityL(0) = \porosityL_0. $$
In case the known initial porosity value is given in Eulerian form, then we consider an initial Jacobian $J_0=1$ such that 
    $$ \porosityL_0 = \porosityE_0.$$
The lack of an initial Jacobian in our formulation stems from the choice of using a quasi-static elasticity. Despite this limitation, we maintain this choice as we found this approach more numerically robust in our preliminary tests due to oscillations arising from the inertia.

We make two remarks on the chosen law of mass conservation. The first one is that we are implicitly assuming that the fluid is incompressible as $\rho_f = \rho_{f0}$. This can be justified by adequate thermodynamic modeling, and is common in the porous media community \cite{chapelle2014} if the pore content is given by an incompressible fluid. The second one is that we are only imposing mass conservation of the fluid phase and not that of the solid phase. This is the case because mass conservation of the solid phase establishes a relationship between the solid density in spatial configuration and deformation, which we do not require in virtue of the quasi-static approximation and the lack of external body forces. Even in the case of inertial deformation, the computation of an Eulerian solid density is typically circumvented by solving the problem in the reference configuration. This would of course not be true for the Eulerian one.

\paragraph{Solid incompressibility.} This condition states that the solid phase is incompressible, and is given in an arbitrary subdomain $\omega_0\subset\Omega_0$ by \cite{Coussy2004}:
    \begin{equation}\label{eq:incompressibility}
    \int_{\omega_0} \left(1-\porosityE_{0}\right)\,dX = \int_\omega \left(1 - \porosityE\right)\,dx,
    \end{equation}
    for some initial porosity $\porosityE_{0}$ with $\omega=\vec x(\omega_0)$. Equation~\eqref{eq:incompressibility} can be rewritten through a localization argument as 
    \begin{equation}\label{eq:incompressibility strong reference}
        \porosityL_s = \porosityE_{s,0} \qquad\text{in $\Omega_0$}, 
    \end{equation}
    where $\porosityL_s$ and $\porosityE_{s,0}$ are the porosities of the solid phase in the Lagrangian and Eulerian form respectively. This yields $J - \porosityL = 1 - \porosityE_{0}$ and can be included in the potential as
    $$ \tilde\Psi \coloneqq \Psi + \lambda\left(J - \porosityL - \left(1 - \porosityE_{0}\right)\right),$$
    where $\lambda$ is the Lagrange multiplier associated to the constraint. As a result, the modified first Piola--Kirchhoff stress tensor is given by
    $$ \widetilde{\ten P} = \parder{\tilde\Psi}{\ten F} = \ten P + \lambda J\ten F^{-T}, $$
    and the pressure becomes
    $$ \widetilde{p} = \parder{\tilde\Psi}{\porosityL} = p - \lambda. $$

This fundamental hypothesis establishes the feedback of the fluid phase on the solid one. To obtain such a feedback without using incompressibility, an additional energy term must be considered in $\Psi$ that is driven by the solid porosity and can be interpreted as the solid phase quasi-incompressibility. See \cite{barnafi2022perfusion} for more details on this topic.
    
\paragraph{Complete forward model.} Putting everything together and assuming homogeneous boundary conditions (either Dirichlet or Neumann), we obtain that the weak formulation of \eqref{eq:momentum}, \eqref{eq:mass}, and \eqref{eq:incompressibility} is given by: Find a displacement $\vec d$ in $V^d$, a Lagrangian porosity $\porosityL$ in $V^\porosityL$, and a Lagrange multiplier $\lambda$ in $V^\lambda$ such that 
    \begin{equation}\label{eq:forward weak}
        \begin{aligned}
            \int_{\Omega_0} \ten P(\ten F, \porosityL): \grad \vec d^*\,dX + \int_{\Omega_0}\lambda J\ten F^{-T}:\grad \vec d^*\,dX &=\int_{\Omega_0}\vec g\cdot \vec d^*\,dX &&\forall \vec d^*\in V^d, \\
            \int_{\Omega_0}\der{\porosityL}{t}\porosityL^*\,dX + \int_{\Omega_0} \ten K(\ten F) \grad \widetilde p(\ten F, \porosityL)\cdot \grad \porosityL^*\,dX &= \int_{\Omega_0}\rho_{f}^{-1}\Theta\porosityL^*\,dX &&\forall \porosityL^*\in V^\porosityL,\\
            \int_{\Omega_0}\left(J - \porosityL\right)\lambda^*\,dX &= \int_{\Omega_0}\left(1 - \porosityE_{0}\right)\lambda^*\,dX && \forall \lambda ^*\in V^\lambda,
        \end{aligned}
    \end{equation}
where the first Piola-Kirchhoff tensor $\ten P$ and the pressure $p$ are computed from the Helmhotz potential $\Psi$ as $\ten P = \parder{\Psi}{\ten F}$ and $p=\parder{\Psi}{\porosityL}$. 


\paragraph{Note: Hidden second order derivatives.}\label{section:second-order-terms}
We note that problem \eqref{eq:mass} has second order order derivatives coming from the pressure gradient, which is given by 
    $$ \grad p(\ten F, \porosityL) = \parder{p}{\ten F}:\grad \ten F + \parder{p}{\porosityL} \grad \porosityL.$$
The conflicting term is $\grad \ten F$ because $H^1$ conforming finite element schemes do not have a well defined second derivative, and in any case the displacement is typically an $\vec H^1$ function, not $\vec H^2$. This difficulty, which we refer to as \emph{primal inconsistency}, will be addressed in Section~\ref{section:mixed} by formulating mixed models through physically significant auxiliary variables.

\subsection{Eulerian field equations}\label{section:inverse}
To define the Eulerian formulation, we follow the inverse displacement model from \cite{govindjee1996computational}. Consider the inverse displacement $\hat{\vec d}$ such that $\vec X(\vec x) = \vec x + \hat{\vec d}(\vec x)$, which yields a strain tensor $\ten f=\ten I + \grad_{\vec x}\hat{\vec d}$ and Jacobian $j = \det \ten f$ that are related to the forward quantities through $\ten f = \ten F^{-1}$ and $j = 1/J$. Interestingly, the elasticity problem in Eulerian configuration yields a problem that completely determines the inverse displacement $\hat{\vec d}$. In what follows, we pose the current configuration equations of linear momentum \eqref{eq:momentum}, mass conservation \eqref{eq:mass}, and solid incompressibility \eqref{eq:incompressibility}. We highlight that in the current configuration, the natural porosity to consider as a variable is the Eulerian one, i.e. $\porosityE$.

\paragraph{Equilibrium equation.} Changing variables in the first Piola--Kirchhoff stress tensor we obtain that $\ten P(\ten F) = \ten P(\ten f^{-1})$, so that the push-forward of the linear momentum equation is given by
    \begin{equation}\label{eq:momentum current}
        -\dive_{\vec x}\left(j\ten P(\ten f^{-1})\ten f^{-1}\right) = \vec g_e \quad \text{in $\Omega$},
    \end{equation}
for a given volumetric load $\vec g_e$, where $\ten \sigma(\ten f)\coloneqq j\ten P(\ten f^{-1})\ten f^{-1}$ is known as the Cauchy stress tensor. Boundary conditions are given as in the Lagrangian case, with the difference that Neumann and Robin boundary conditions are written using spatial quantities as $\ten\sigma\vec n=\ten t$.

\paragraph{Mass conservation.} We recast mass conservation in terms of its natural Eulerian variable $\porosityE=\porosityL/J = j\porosityL$, and denote the Eulerian fluid velocity $\vec u = J\ten F^{-T}\vec U$ to write mass conservation as
    \begin{equation}\label{eq:mass current}
        \der{\porosityE}{t} + \dive_{\vec x}\vec u = \frac{1}{\rho_{f}}\theta\quad\text{in $\Omega$},
    \end{equation}
for a given source/sink term $\theta$. See \cite{Coussy2004} for a derivation of this equation. As before, the fluid velocity is computed from Darcy's law
    $$ \vec u \coloneqq -\ten k \grad_{\vec x}p.$$
Boundary conditions are given as in the Lagrangian case. The choice of the initial conditions will be discussed in Section~\ref{section:refconf-model} as it is one of the main contributions of this work.

\paragraph{Solid incompressibility.} Pushing forward \eqref{eq:incompressibility} to $\Omega$, we obtain the following Eulerian form of solid incompressibility:
    \begin{equation}\label{eq:incompressibility strong current}
    1 - \porosityE = j\left( 1 - \porosityE_{0}\right)\qquad \text{in $\Omega$}.
    \end{equation}
This is incorporated in our formulation by modifying the Helmholtz potential as
    $$ \tilde \Psi = \Psi + \lambda\left( 1 - \porosityE - j\left(1 - \porosityE_{0}\right)\right), $$
where $\lambda$ is the Lagrange multiplier associated to the constraint. With it, the Cauchy stress tensor becomes
    $$ \widetilde{\ten \sigma} = \ten \sigma + \lambda\ten I, $$
and the pressure remains the same. Note that the stress generated by the Lagrange multiplier is significantly simpler in the Eulerian formulation.

\subsection{Constitutive modeling}\label{section:constitutive}
Even though our formulation is valid independently of the mechanical material model used, we will consider a material related to cardiac modeling, since cardiac applications are our main interest, and are also numerically challenging. In virtue of this, the Helmholtz potential $\Psi$, as in \cite{barnafi2022perfusion}, will be separated additively into a mechanical contribution $\Psi_M$ and porous contributions $\Psi_P$ as
    $$ \Psi(\ten F, \porosityL) = \Psi_M(\ten F) + \Psi_P(\porosityL), $$
where the term $\Psi_P$ acts as a barrier function such that $\porosityL>0$. Similarly to \cite{chapelle2014}, the inequality $\porosityE<1$ (or equivalently $\porosityL<J$) is naturally satisfied using solid incompressibility. Indeed, using equation \eqref{eq:incompressibility strong current} 
    $$\porosityE = 1- j\porosityE_{s,0} \qquad\quad\text{(equivalently $\porosityL = J - \porosityE_{s,0}$)}, $$
and assuming an arbitrarily small, but positive, constant $\epsilon$ such that $\porosityE = 1-\epsilon$, we obtain 
    $$ \epsilon = j\porosityE_{s,0} = J^{-1}\porosityE_{s,0}. $$
Thus, $\epsilon \to 0$ if $j\to 0$, i.e. when $J\to \infty$, which is nonphysical and already characterized by the potential $\Psi_M$. In this work, we consider an Usyk energy of ventricular muscle \cite{Usyk2002} for $\Psi_M$, given by
    $$ \Psi_M(\ten F) =  C\left(\exp(Q(\overline{\ten E}))-1\right) + \frac B 2 (J-1)\log J, $$
where
\begin{equation*}
    \begin{aligned}
    &Q(\overline{\ten E}) =
      b_{\mathrm{ff}} \overline{E}_{\mathrm{ff}}^2  + b_{\mathrm{ss}} \overline{E}_{\mathrm{ss}}^2 + b_{\mathrm{nn}} \overline{E}_{\mathrm{nn}}^2+ b_{\mathrm{fs}} \left( \overline{E}_{\mathrm{fs}}^2 + \overline{E}_{\mathrm{sf}}^2 \right) + b_{\mathrm{fn}} \left( \overline{E}_{\mathrm{fn}}^2 + \overline{E}_{\mathrm{nf}}^2 \right) + b_{\mathrm{sn}} \left( \overline{E}_{\mathrm{sn}}^2 + \overline{E}_{\mathrm{ns}}^2 \right),\\
    &\overline{E}_\text{ab} =
    \overline{\ten E} \vec {a} \cdot \vec {b}, \qquad a, b \in \{ f, s, n \}, \qquad \overline E = \frac 1 2\left( \overline{\ten F}^T \overline{\ten F} - \ten I\right), \qquad \overline{\ten F} = J^{1/3}\ten F,
    \end{aligned}
\end{equation*}
where  $(\vec f, \vec s, \vec n)$ is a locally orthonormal basis that represents the orientation of the cardiac fibers, $C=880\,\texttt{Pa}$ and $B=5\times 10^4\,\texttt{Pa}$. For the porous energy $\Psi_P$ we consider the arterial Bruinsma et al. energy \cite{bruinsma1988model}, upscaled in \cite{cookson2012novel} to be given by
    $$ \Psi_P(\porosityL) = \frac{q_{1}}{q_{3}}\exp \left(q_{3}\porosityL\right)+q_{2}\porosityL\left(\log\left(q_{3}\porosityL\right)-1\right) $$
where the values of the parameters $q_{i}$ for arteries, capillaries, and veins have been computed in \cite{bruinsma1988model}. Still, as observed in \cite{barnafi2022perfusion}, the energy $\Psi_P$ must be normalized such that it is minimized at the reference values of $\porosityL$, i.e. the pressure at rest yields the reference pressures. Therefore, we make the dependence of the energy explicit on the reference porosities as
    \begin{equation}\label{eq:normalized energy}
        \widetilde{\Psi}_P(\porosityL, \porosityL_{0}) \coloneqq \Psi_P(\porosityL) - \parder{\Psi_P}{\porosityL}(\ten I, \porosityL_0) \porosityL.
    \end{equation}
By considering also the solid incompressibility constraint, we obtain the following expressions for the first Piola--Kirchhoff stress tensor and the pressure:
    \begin{align*}
        \ten P(\ten F) &= \parder{\Psi_M}{\ten F}(\ten F) + \lambda J \ten F^{-T}, \\
        p(\porosityL, \porosityL_{0}) &= \parder{\Psi_P}{\porosityL}(\porosityL) - \parder{\Psi_P}{\porosityL}(\porosityL_0) - \lambda.
    \end{align*}
We make the following very important remarks.
\begin{remark}
    The first Piola--Kirchhoff tensor and the pressure seem decoupled from one another. Indeed, it is the Lagrange multiplier that yields the interaction between fluid and solid phases, and it can be interpreted as the pressure existing at the interface between both phases. 
\end{remark}
\begin{remark}
    The pressure contribution from the reference porosity $\porosityL_{0}$ has a flipped sign with respect to the pressure obtained from the energy. This flipped sign, together with the negative sign obtained from the time-reversal strategy proposed, are what yield that the nonlinear porous media equation behaves as a heat equation instead of a reversed heat equation.
\end{remark}
\begin{remark}
    The pressure is considered such that it is gives a reference porosity of 0 at $\porosityL_0$ for simplicity. It could of course be considered that the reference porosity is non-zero, i.e.
        $$ p(\porosityL, \porosityL_0) = \parder{\Psi_P}{\porosityL}(\porosityL) - \parder{\Psi_P}{\porosityL}(\porosityL_0) + p_\texttt{ref} - \lambda, $$
    for some reference pressure $p_\texttt{ref}$.
\end{remark}

\subsection{Formulation of the reference configuration model}\label{section:refconf-model}
In this section we extend the elastic reference configuration model to poroelasticity. The extension is based on two main points: (i) the notion of mechanical equilibrium can be understood as a steady state solution in the context of porous media, and (ii) the resulting steady state equation can present extreme nonlinearities that are very challenging to solve computationally. This motivates considering another related time dependent problem instead of the steady state problem directly. We now provide the details of each point separately.

\paragraph{Equilibrium.} The reference configuration problem departs from the assumption that we know not only the forces acting of a system, but also its \emph{resulting} solution. The time dependent nature of system \eqref{eq:mass current} makes it unclear at exactly which time we observe the system, so we assume that the solution given is \emph{stationary}, i.e. it is the solution $(\hat{\vec d}, \porosityE_0)$ of \eqref{eq:momentum current} and the following stationary mass conservation problem:
    \begin{equation}\label{eq:equilibrium}
        -\dive_{\vec x} \ten k\grad_{\vec x}p(\ten f^{-1}, j\porosityE, \porosityE_0) = \frac 1{\rho_f} \theta \qquad\text{on}\quad\Omega,
    \end{equation}
as seen from \eqref{eq:mass current}. We highlight that in this model the given porosity is $\porosityE$, and the unknown is the reference one $\porosityE_0$. The way in which $\theta$ is considered here is problem dependent. For example, if it varies in time in the forward model, then its initial value could be considered. If it is instead an autonomous function of the variables (as we do in the numerical tests section), i.e. $\theta = \theta(\ten f^{-1}, j\porosityE, \porosityE_0)$, then it can be left as it is. 

\paragraph{Time dependence.} As seen from Section~\ref{section:constitutive}, the pressure function can be severely nonlinear, and indeed all of our preliminary attempts at solving problem \eqref{eq:equilibrium} resulted in diverging nonlinear iterations. For this reason, and to allow for more flexible continuation methods, we provide a way to reconsider time dependence in \eqref{eq:equilibrium} in a way that yields the reference configuration solution. Equation \eqref{eq:mass current} can be rewritten as 
    $$ \parder{(\porosityE - \porosityE_0)}{t} - \dive_{\vec x} \ten k \grad_{\vec x} p = \frac 1{\rho_f}\theta, $$
as $\porosityE_0$ is constant in time and the quantity $\porosityE - \porosityE_0$ can be regarded as the variation of porosity. Inverting the variables, we consider a given (time constant) porosity $\porosityE=\overline{\porosityE}$ and write the reference mass conservation problem as
    \begin{equation}\label{eq:mass backwards}
    \begin{aligned}
        - \parder{\porosityE_0}{t} - \dive_{\vec x} \ten k \grad_{\vec x} p(\ten f^{-1}, j\overline \porosityE, \porosityE_0) &= \frac 1{\rho_f}\theta &&\text{on}\quad\Omega, \\
        \porosityE_0(0) &= \overline \porosityE &&\text{on}\quad\Omega.
    \end{aligned}
    \end{equation}

Naturally, a steady-state solution of problem \eqref{eq:mass backwards} provides a solution to \eqref{eq:equilibrium}.  If we think of \eqref{eq:mass backwards} as a heat equation for a moment, the reference mass conservation equation behaves as a reversed heat equation, which is severely ill-posed. This is not the case here because the normalization factor $\parder{\Psi_P}{\porosityL}(\porosityL_0)=\parder{\Psi_P}{\porosityL}(\porosityE_0)$ shown in Section~\ref{section:constitutive} has a negative sign.

\paragraph{Complete reference configuration model.} Putting everything together, we obtain the weak formulation of \eqref{eq:momentum current}, \eqref{eq:mass backwards} and \eqref{eq:incompressibility strong current} assuming homogeneous Dirichlet or Neumann boundary conditions as: Find an inverse displacement $\hat{\vec d}$ in $V^d$, a spatial reference porosity $\porosityE_0$ in $V^\porosityE$, and a Lagrange multiplier $\lambda$ in $V^\lambda$ such that
    \begin{equation}\label{eq:backward weak}
        \begin{aligned}
            \int_\Omega \left(j\ten P(\ten f^{-1})\ten f^{-T}\right):\grad_{\vec x}\vec d^*\,dx + \int_\Omega \lambda \dive_{\vec x} \vec d^*\,dx &=\int_\Omega \vec g_e\cdot \vec d^*\,dx &&\forall \vec d^* \in V^d, \\
            -\int_\Omega \der{\porosityE_0}{t}\porosityE^*\,dx + \int_\Omega \ten k \grad_{\vec x} \widetilde p(\ten f^{-1}, j\overline\porosityE, \porosityE_0):\grad_{\vec x}\phi^*\,dx &= \int_\Omega \rho_f^{-1}\theta\porosityE^*\,dx && \forall \porosityE^*\in V^\porosityE, \\
            \int_\Omega j\left(1 - \porosityE_0\right)\lambda^* \,dx &= \int_\Omega \left(1 - \overline\porosityE\right)\lambda^*\,dx && \forall \lambda^*\in V^\lambda, 
        \end{aligned}
    \end{equation}
where $\ten P = \parder{\Psi}{\ten F}$, $p = \parder{\Psi}{\porosityL}$, and $\porosityE_0(0)=\overline\porosityE$. 

\section{Numerical solution strategy}\label{section:numerical}
In this section we provide numerical details for solving the reference configuration problem. All computations were performed using the Firedrake library \cite{firedrake2016}. We briefly provide the discretization choices used to solve both \eqref{eq:forward weak} and \eqref{eq:backward weak}:

\paragraph{Space discretization.} We use the inf-sup stable finite element space $\mathbb{P}_2\times \mathbb P_1 \times \mathbb P_1$ \cite{barnafi2022edema}, where higher order elements are used for the displacement. Firedrake ships Gauss-Jacobi quadrature rules by default, which we have set to a degree of 6 to handle the nonlinearities adequately. Anything lower than that exhibits numerical artifacts in our experience.

\paragraph{Time discretization.} Given the quasi-static nature of most equations, we only need to provide a time-discretization for the porous media equations. These equations are parabolic (in virtue of the convex potential) \cite{vazquez2007porous}, so we use a simple backward Euler scheme \cite{hairer1993solving}, and thus treat all variables implicitly.

\paragraph{Iterative solvers.} The nonlinear system arising from the implicit discretization is solved using Newton's method \cite{wright1999numerical}, where the linear system is solved at each time step using the direct solver available in the MUMPS library \cite{amestoy2000mumps}, which includes an MPI-parallel implementation. Better iterative solvers for nonlinear poroelasticity based on Krylov Subspace methods can be found in \cite{barnafi2022coupling}, but their extension to the fully nonlinear case we considered has not yet been developed.

\paragraph{Continuation.} The model considered is highly nonlinear, and thus a naive use of Newton's method yields diverging iterations. To avoid this, we consider a continuation strategy for the source term, such that for a given source term $\theta^*$ (possibly time dependent) we consider the linear ramp given by $\theta(t) = \min\{t/t_\text{ramp}, 1\}\theta^*(t)$, so that the source term is fully active at $t=t_\text{ramp}$.

\paragraph{Convergence criterion.} To compute the reference solutions as stationary states of equation \eqref{eq:backward weak}, we solve the equations until the normalized residual of the steady state mass conservation equation \eqref{eq:equilibrium} is below a given threshold. We thus set
    $$ R_i(t) \coloneqq \int_\Omega \ten k\grad_{\vec x} p\cdot \grad_{\vec x} \phi_i^* - \rho_f^{-1}\theta(t)\phi_i^*\,dx, $$
where $\phi_i^*$ stands for the $i$-th test function in the discrete space, and use as a convergence criterion the following:
    $$  \left| \vec R(t)\right|_{\ell^2} \leq \texttt{tol}R_0, $$
where $[\vec R]_i = R_i$, $|(\cdot)|_{\ell^2}\coloneqq \sqrt{\sum_i(\cdot)^2_i}$, and $R_0\coloneqq |R(t_\text{ramp})|_{\ell^2}$ is a normalization factor. We assume a tolerance of $\texttt{tol}=10^{-6}$ was used unless otherwise stated.

\subsection{Accelerated stationary state computation}\label{section:acceleration}
In this section we show that the computation of a stationary state can be accelerated through a fixed-point acceleration technique known as Anderson acceleration \cite{anderson1965iterative}. The method changes a fixed point iteration $x^k = g(x^{k-1})$ with a memory dependent algorithm given by 
    $$ x^k =\sum_{i\in\{1,..., m\}} \alpha_i g(x^{k-i}), $$
where the $\alpha_i$ weights are optimal for a given norm, and $m$ is known as the depth parameter. More details on how to compute the weights, as well as the theoretical foundations of this technique, can be found in \cite{walker2011anderson}. We highlight that this method is equivalent to the GMRES method in the linear case, and to the multisecant method in the nonlinear case.

To apply Anderson acceleration to the problem of finding a stationary state, we formulate a fixed-point operator in two steps. First, we denote with $\mathbb T: V^i\to V^d\times V^\porosityE \times V^\lambda$ the operator given by 
        $$ (\hat{\vec d}^k, \porosityE^k, \lambda^k) = \mathbb T(\porosityE^{k-1}), $$
where $(\cdot)^k$ is the discrete solution at instant $t^k = k \Delta t$. Then we denote with $\Pi_\porosityE$ the projection 
        $$ \porosityE = \Pi_\porosityE(\hat{\vec d}, \porosityE, \lambda). $$
With these operators, we note that the problem of finding a stationary state to equation \eqref{eq:backward weak} can be recast as finding a fixed-point of the operator $\mathbb S: V^\porosityE\to V^\porosityE$ given by
        $$\mathbb S = \Pi_\porosityE \circ \mathbb T. $$
This operator is amenable to acceleration, and we compare the number of iterations it takes to compute the stationary state with and without Anderson acceleration in Section~\ref{section:tests}. Naturally, all computations performed here can be trivially extended to the forward problem. 

\section{Consistent mixed models}\label{section:mixed}
In this section we propose two mixed models associated with nonlinear poroelasticity that can circumvent the primal inconsistency problem, detailed in Section~\ref{section:second-order-terms}. The first one is the commonly used mixed formulation given by considering the fluid velocity as a variable $\vec u\coloneqq -\ten \kappa \grad p$, which is useful also for extended models in which the Brinkman effect is included in the fluid inertia. The second one is a simpler formulation in which the pressure becomes a variable $\mu\coloneqq p(\ten F, \porosityL)$. Given that the main difficulty is the treatment of the mass conservation equation, throughout this section we focus on the following equation: 
    \begin{equation}\label{eq:porous media}
        \der{\porosityE}{t} - \dive \left(\ten k\grad p(\ten F, \porosityE)\right) = \theta \qquad\text{in $\Omega$},
    \end{equation}
for given functions $\ten F$ and $\theta$. We highlight that this problem is independent of the previous models, so we consider $\Omega$ to be a general domain in this section, and we avoid unnecessary notation regarding the differentiation variables. The boundary conditions for each formulation will be different, thus, to highlight how to convert from one formulation to the other one, we consider mixed boundary conditions on Dirichlet and Neumann boundaries, $\Gamma_D$ and $\Gamma_N$ respectively, with
    $$ \porosityE = \porosityE_D \quad\text{on $\Gamma_D$}, \quad \text{and} \quad \ten k \grad p\cdot \vec n = \bar p \quad\text{on $\Gamma_N$}, $$
where it holds that $\overline{\partial\Omega} = \overline{\Gamma_D}\cup \overline{\Gamma_N}$.

\subsection{Mixed formulation in pressure}
Starting from equation \eqref{eq:porous media}, we define the variable $\mu \coloneqq p(\ten F, \porosityE) $ so that we can rewrite the problem as the following system of equations: 
    \begin{equation}\label{eq:mixed-p strong}
        \begin{aligned}
        \der{\porosityE}{t} - \dive \ten k \grad \mu &= \theta &&\text{in $\Omega$},\\
        \mu - p(\ten F, \porosityE) &=0 &&\text{in $\Omega$}.
        \end{aligned}
    \end{equation}
From the integration by parts of the second order operator against a test function $\mu^*$:
    $$ -\int_\Omega \dive \left(\ten k \grad \mu\right)\mu^*\,dx = \int_\Omega \ten k \grad\mu\cdot\grad\mu^*\,dx - \int_{\partial\Omega} \left(\ten k \grad\mu\cdot \vec n\right)\mu^*\,ds $$
    we can see that the Dirichlet boundary conditions for this model are given by 
    $$ \mu = \mu_D \quad\text{on $\Gamma_D$}, $$
    which in virtue of the definition of $\mu$, correspond to imposing pressure values on the boundary instead of porosity values as in the primal formulation \eqref{eq:porous media}. The Neumann boundary condition is instead given by 
    $$ \ten k\grad \mu\cdot \vec n = \bar p \quad\text{on $\Gamma_N$}, $$
    which in virtue of the definition of $\mu$, correspond to the same Neumann boundary conditions used in \eqref{eq:porous media}, thus justifying the use of the same function $\bar p$. We highlight that these boundary conditions are better from a physical point of view, as typically in real applications the boundary values correspond to pressure instead of porosity. 

    The weak formulation of \eqref{eq:mixed-p strong} is given by finding $\porosityE$ in $V^\porosityE$, and $\mu$ in $V^\mu$ such that 
    \begin{equation}\label{eq:mixed-p weak}
        \begin{aligned}
            \int_\Omega \der{\porosityE}{t}\mu^*\,dx + \int_\Omega \ten k\grad \mu\cdot \grad\mu^*\,dx &= \int_\Omega\theta \mu^*\,dx &&\forall \mu^*\in V^\mu, \\
            \int_\Omega \mu \porosityE^*\,dx - \int_\Omega p(\ten F, \porosityE)\porosityE^*\,dx &= 0 &&\forall \porosityE^* \in V^\porosityE.
        \end{aligned}
    \end{equation}
    This formulation does not present a gradient on the nonlinear function $p$, thus is avoids the primal inconsistency problem.

\begin{remark}
The choice of switching the test functions is purely theoretical. To see this, consider the simplest case where $p(\ten F, \porosityE) = \porosityE$ and an implicit time discretization with time-step $\Delta t = 1$. If we use the test functions $(\mu^*, \porosityE^*)=(\mu, -\porosityE)$ in \eqref{eq:mixed-p weak} and add both equations, we obtain
    $$\int_\Omega |\grad \mu |^2\,dx + \int_\Omega \porosityE^2\,dx, $$
which hints on better theoretical properties for this formulation, as well as not requiring an inf-sup condition among the spaces $V^\porosityE$ and $V^\mu$. This in particular allows for a robust lowest order discretization.
\end{remark}

\subsection{Mixed formulation in fluid velocity}
Proceeding as before, we define the function $\vec u \coloneqq - \ten k \grad p$, which converts \eqref{eq:porous media} into the following system of equations:
    \begin{equation}\label{eq:mixed-u strong}
        \begin{aligned}
            \der{\porosityE}{t} + \dive \vec u &= \theta &&\text{in $\Omega$}, \\
            \ten k^{-1}\vec u + \grad p &= 0 &&\text{in $\Omega$}.
        \end{aligned}
    \end{equation}

Again, we deduce the boundary conditions from the integration by parts of the $\grad p$ term against a test function $\vec u^*$: 
    $$ \int_\Omega \vec u^*\cdot \grad p\,dx = -\int_\Omega p\dive\vec u^*\,dx + \int_{\partial\Omega} p \vec u^*\cdot \vec n\,dx.$$
In this case, the Dirichlet boundary conditions correspond to 
    $$ \vec u\cdot \vec n = -\bar p \quad\text{on $\Gamma_D$}, $$
and the Neumann ones correspond to
    $$ p = p_N \quad\text{ on $\Gamma_N$}.$$
Interestingly, the boundary conditions of this formulation correspond to the ones from \eqref{eq:mixed-p strong}, but with exchanged roles (i.e. Dirichlet conditions become Neumann conditions and vice-versa). This phenomenon is typical in mixed methods \cite{gatica2014simple}, but the nonlinearity of $p$ breaks that relationship with respect to the original primal formulation \eqref{eq:porous media}. 

Its weak formulation is given by finding $\vec u$ in $V^u$ and $\porosityE$ in $V^\porosityE$ such that 
    \begin{equation}\label{eq:mixed-u weak}
        \begin{aligned}
            \int_\Omega \der{\porosityE}{t}\porosityE^*\,dx + \int_\Omega \porosityE^*\dive \vec u\,dx &= \int_\Omega \theta\porosityE^*\,dx  && \forall \porosityE^*\in V^\porosityE, \\
            \int_\Omega \ten k^{-1}\vec u\cdot \vec u^*\,dx - \int_\Omega p(\ten F, \porosityE) \dive \vec u^*\,dx &= 0 &&\forall \vec u^*\in V^u.
        \end{aligned}
    \end{equation}
Similarly to \eqref{eq:mixed-p weak}, this formulation does not present a gradient on the nonlinear function $p$, and thus it also avoids the primal inconsistency error.

\begin{remark}
The mixed $\vec u$ formulation presents two main advantages: (i) it yields the physically significant variable $\vec u$ with optimal convergence rates, and (ii) it allows for generalized Darcy models such as the ones considering a Brinkman term, or also fluid inertia. Still, this comes with the disadvantage of being a more difficult problem. Indeed, considering the test functions $(\vec u^*, \porosityE^*) = (\vec u, \porosityE)$ as before, we obtain 
    $$ \int_\Omega |\vec u|^2\,dx + \int_\Omega |\porosityE|^2\,dx, $$
but the presence of the divergence operator shows we require $\vec u$ to be in $H(\dive; \Omega)$. This can be overcome using the theory of saddle point problems, which also shows that inf-sup stable finite element spaces need to be chosen for the $V^u\times V^\porosityE$ pair. 
\end{remark}

\begin{remark}
The mixed models follow the same logic used to define the fixed-point operator in Section~\ref{section:acceleration}, meaning that they are also amenable to acceleration.
\end{remark}
We show the main features of all three models in Table~\ref{table:models}. We emphasize the difference in terms of the variables and the boundary conditions, where again we highlight the following: (i) The mixed models present more physically adequate boundary conditions, as the quantity that is typically measured is the pressure, not the porosity, (ii) the role of the boundary conditions between the mixed models is exchanged, and (iii) the mixed models do not require second order derivatives of the displacement.

\begin{table}[ht!]
    \centering
    \begin{tabular}{r | c c c c }
    \toprule 
    Model     & Variables              & Dirichlet BC                    & Neumann BC                            & Consistent \\\midrule
    Primal    & $\porosityE$           & $\porosityE=\bar\phi$           & $\grad p \cdot \vec n = \bar{\vec u}$ & \xmark \\
    Mixed-$p$ & $\porosityE, \mu$      & $\mu = \bar p$                  & $\grad p\cdot \vec n = \bar{\vec u}$  & \cmark \\
    Mixed-$u$ & $\porosityE, \vec u$   & $\vec u\cdot\vec n=\bar{\vec u}$& $p = \bar p$                          & \cmark\\\bottomrule 
    \end{tabular}
    \caption{Comparison of variables and boundary conditions required for the primal, pressure mixed and velocity mixed models. }
    \label{table:models}
\end{table}

\section{Numerical tests}\label{section:tests}
In this section we provide numerical tests that support our theory. The models have been implemented in Firedrake \cite{firedrake2016} following the numerical solution strategy detailed in Section~\ref{section:numerical}, and the visualizations are performed with Paraview \cite{paraview2015}. The tests considered, and their scope, are as follows:
    \begin{description}
        \item[Test 1:] We compute the reference configuration of a 2D geometry with a physically motivated nonlinear source term. The scope of this test is to show that the loaded reference configuration yields the given current configuration, and also recovers the given porosity.
        \item[Test 2:] We repeat Test 1 for a 3D geometry. 
        \item[Test 3:] We verify the impact of using Anderson acceleration to compute the desired stationary state in both 2D and 3D tests, so we perform acceleration for varying levels of memory and assess its performance by measuring the number of iterations required for convergence.
        \item[Test 4:] We verify that the mixed models are also able to yield a satisfactory solution to the reference configuration problem by comparing the solution obtained by all 3 models (primal, pressure mixed, and velocity mixed).
        \item[Test 5:] We compute the reference configuration of a realistic geometry using physically accurate parameters. 
    \end{description}

In both square and slab tests, we use the canonical basis as fiber directions. Instead, for the left ventricle simulation we use physically accurate fiber directions, computed numerically as in \cite{barnafi2023modeling}. More details on this are provided in the corresponding test.

\subsection{Test 1: Validation in 2D}
In this test, we validate the performance of the reference configuration formulation. The way to do this is by loading the problem with a pressure driven source term that models the inflow of fluid according to a pressure difference:
    $$ \theta = -\beta(p(\ten F, \porosityL) - p_a), $$
where $p_a=10^4\,\texttt{Pa}$ and $\beta=10^{-4}\,\frac{1}{\texttt{s Pa}}$. For the porous media potential we used the arterial pressure model described in Section~\ref{section:constitutive} where $q_1 = 1.333\,\texttt{Pa}$, $q_2=550\,\texttt{Pa}$, and $q_3 = 10$, with a given porosity of $\overline \phi=0.1$. Finally, we considered an isotropic permeability given by the scalar $k=2\cdot\,10^{-7}\texttt{m}^2 \texttt{(s\,Pa)}^{-1}$.  We considered a square geometry $\Omega$ of side length $1\,\texttt{cm}$, a time step of $\Delta t=0.01\,s$, and a ramp time given by $t_\text{ramp}=0.1\,s$. As boundary conditions, we used Dirichlet conditions that allowed for sliding, i.e. $\vec d_x = 0$ on $\{x=0\}$ and $\vec d_y = 0$ on $\{y=0\}$, and homogeneous Neumann conditions elsewhere. 

\begin{figure}[ht!]
    \centering
    \begin{subfigure}{0.49\textwidth}
        \includegraphics[width=\textwidth]{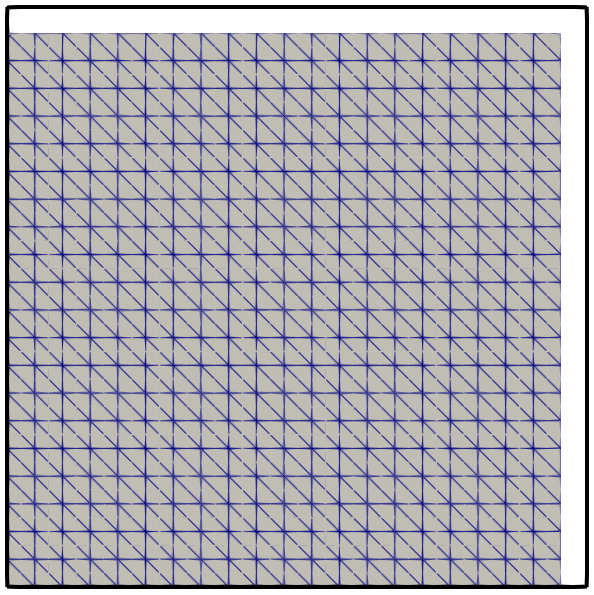}
        \caption{}
    \end{subfigure}
    \begin{subfigure}{0.49\textwidth}
        \includegraphics[width=\textwidth]{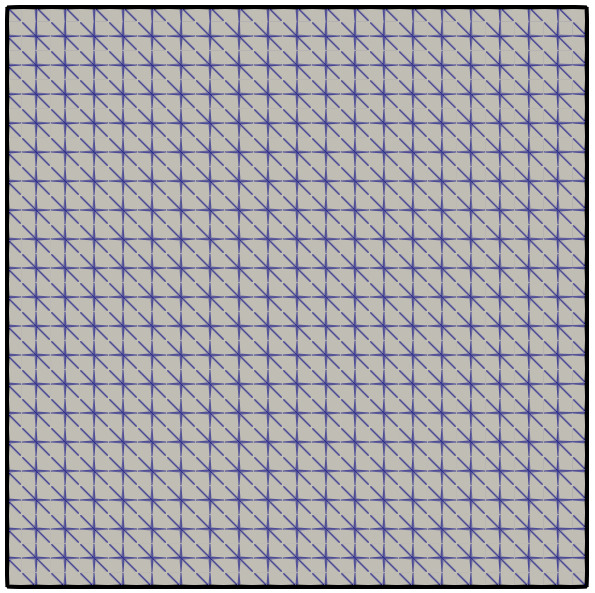}
        \caption{}
    \end{subfigure}
    \caption{(a) the computed reference configuration with respect to the initial geometry (black contour). (b) the forward problem solution as computed starting from the reference configuration computed in (a). }
    \label{fig:test-one-phase}
\end{figure}

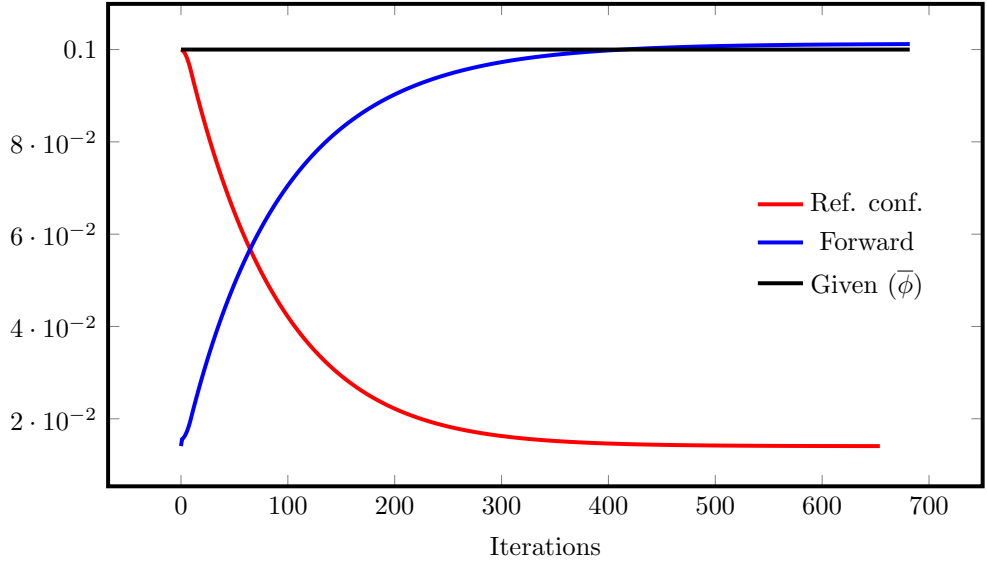
\begin{figure}[ht!]
    \centering
    \begin{tikzpicture}
        \begin{axis}
            [width=0.8\textwidth, height=8cm, xlabel=Iterations, ylabel=, tick label style={font=\normalsize},line width=1.5pt,legend style={draw=none, fill opacity=0.7, text opacity= 1,row sep=2pt, font=\normalsize, at={(0.95,0.5)}, anchor=east}, legend columns=1]
        \addplot+[mark=none,color=red] table [x=Time, y=phiAvg, col sep=comma] {pm-one-phase-back.csv};
        \addplot+[mark=none,color=blue] table [x=Time, y=phiAvg, col sep=comma] {pm-one-phase-forw.csv};
        \addplot+[mark=none,color=black] coordinates {(0,0.1) (682, 0.1)};
        \legend{Ref. conf., Forward, Given ($\overline \phi$)}
        \end{axis}
    \end{tikzpicture}
    \caption{Evolution of the average spatial porosity $\frac 1{|\Omega|}\int_\Omega\phi\,dx$ along the iterations of both the reference configuration and forward problems. The given value $\overline \phi$ is shown for comparison.}
    \label{fig:test-one-phase-avg}
\end{figure}

In Figure~\ref{fig:test-one-phase}, we show the deformed solution in two snapshots: In the left figure, we show the computed reference configuration with respect to the initial geometry. In the right one, we show the solution obtained by solving the forward problem on the computed reference configuration, where the accuracy of this method can be appreciated. Additionally, we computed the average porosity $\frac 1{|\Omega|}\int_\Omega\phi\,dx$ and displayed its evolution on the current configuration in Figure~\ref{fig:test-one-phase-avg}, where we plotted the porosity for both reference configuration and forward problems. Here, it is possible to observe that the reference configuration deflates and loses a significant amount of fluid, given by the 84\% of the initial value as it arrives at an average of 0.014. The recovered porosity presents a small error of roughly a 1.1\% with respect to the correct value, which we attribute to approximation error.

\subsection{Test 2: Validation in 3D}

\begin{figure}[ht!]
    \centering
    \begin{subfigure}{0.49\textwidth}
        \includegraphics[width=\textwidth]{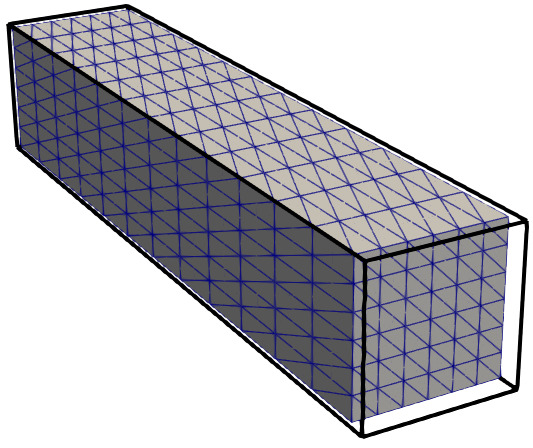}
        \caption{}
    \end{subfigure}
    \begin{subfigure}{0.49\textwidth}
        \includegraphics[width=\textwidth]{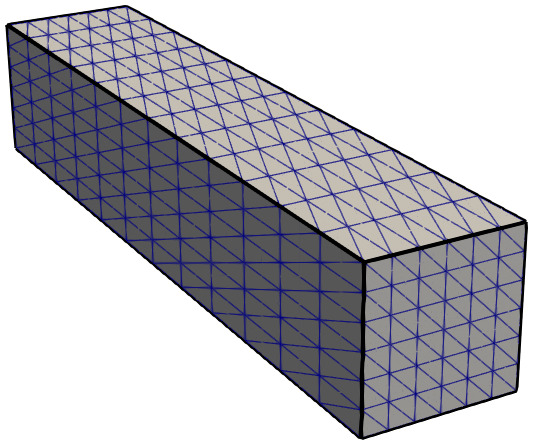}
        \caption{}
    \end{subfigure}
    \caption{(a) the computed reference configuration with respect to the initial geometry (black contour). (b) the forward problem solution as computed starting from the reference configuration in (a). }
    \label{fig:test-3d}
\end{figure}

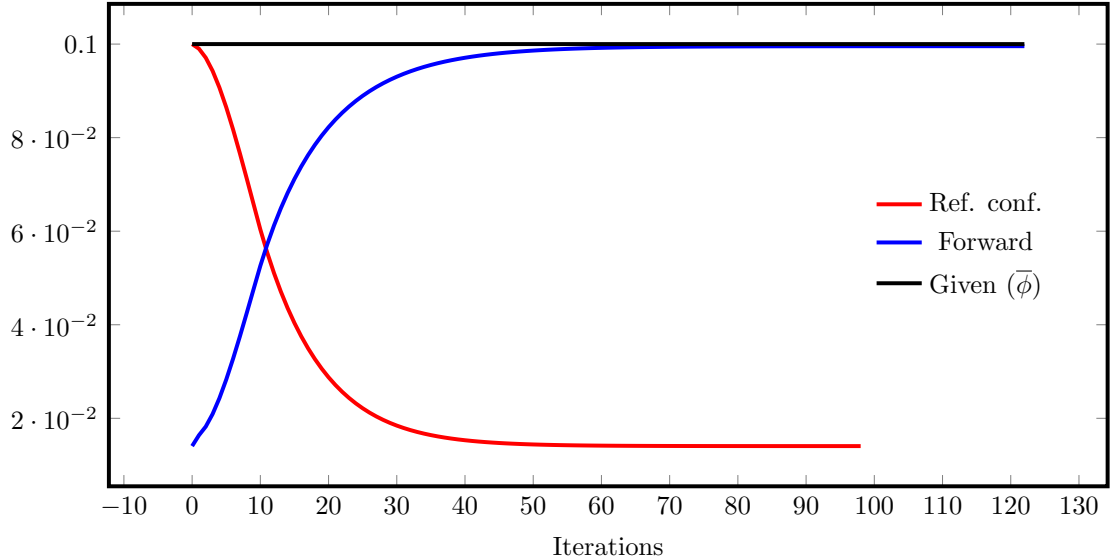
\begin{figure}[ht!]
    \centering
    \begin{subfigure}{0.9\textwidth}
    \begin{tikzpicture}
        \begin{axis}
            [width=\textwidth, height=8cm, xlabel=Iterations, ylabel=, tick label style={font=\normalsize},line width=1.5pt,legend style={draw=none, fill opacity=0.7, text opacity= 1,row sep=2pt, font=\normalsize, at={(0.95,0.5)}, anchor=east}, legend columns=1]
        \addplot+[mark=none, color=red] table [x=Time, y=phiAvg, col sep=comma] {pm-one-phase-3d-back.csv};
        \addplot+[mark=none, color=blue] table [x=Time, y=phiAvg, col sep=comma] {pm-one-phase-3d-forw.csv};
        \addplot+[mark=none, color=black] coordinates {(0,0.1) (122, 0.1)};
        \legend{Ref. conf., Forward, Given ($\overline \phi$)}
        \end{axis}
    \end{tikzpicture}
    \caption{Phase 1}
    \end{subfigure}
    \caption{Evolution of the average spatial porosity $\frac 1{|\Omega|}\int_\Omega\phi\,dx$ along the iterations of both the reference configuration and forward problems in both phases for the brick 3D geometry. The given value $\overline \phi$ is shown for comparison.}
    \label{fig:test-3d-avg}
\end{figure}

Similar to test 1, we compute the reference configuration for a 3D slab geometry given by $\Omega=(0,5\,\texttt{cm})\times (0,1\,\texttt{cm})\times (0,1\,\texttt{cm})$. The only difference with respect to the 2D test is the relative tolerance $\texttt{tol}=10^{-5}$. In Figure~\ref{fig:test-3d}, we show the deformed solution in two snapshots as before: reference configuration and forward problems in the left and right figures respectively. As in the previous case, the forward solution presents an excellent match with respect to the initial geometry. The reference porosities are displayed in Figure~\ref{fig:test-3d-avg}, where the results obtained are analogous to the ones from Test 1.

\subsection{Test 3: Acceleration test} 
In this section, we test the impact of using Anderson acceleration to reduce the iterations required for convergence in both Tests~1 and 2. The number of iterations required for both the reference configuration and forward problems is shown in Table~\ref{table:acceleration}, with respect to different levels of acceleration. For the sake of comparison, we used a tolerance of $\texttt{tol}=10^{-5}$ for both 2D and 3D tests.

Acceleration is not used immediately, but only when the ramping has finished in order to accelerate with respect to the desired dynamics, which happens after 10 timesteps. As seen from Table~\ref{table:acceleration}, we make the following observations: (i) both formulations require roughly 100 iterations for convergence without acceleration, (ii) acceleration works even using only 1 previous vector, and (iii) all accelerated formulations take roughly less than 80\% the number of iterations it takes to solve the problem, where using only one vector yields the lowest iteration count. We conclude that acceleration is an effective strategy to reduce the overall computational cost of computing the reference configuration as a stationary state of the nonlinear porous media equations.

    \begin{table}[ht!]
        \centering
        \begin{tabular}{c | c c c c}
        \toprule &  AA(0) &  AA(1)  &  AA(2)  &  AA(5)  \\\midrule
        Ref. conf. 2D&  98  &  15     &  17     &  20     \\
        Forward  2D&  122 &  14     &  15     &  18     \\\midrule
        Ref. conf. 3D&  98  &  15     &  17     &  20     \\
        Forward  3D&  122 &  14     &  15     &  18     \\\midrule
        \end{tabular}
        \caption{Test 3, number of time steps required to converge to a stationary state using Anderson Acceleration for both forward and reference configuration models.}
        \label{table:acceleration}
    \end{table}

\subsection{Test 4: Validation of mixed models} 
\begin{figure}[ht!]
    \centering
    \newcommand{\hei}{6cm}
    \begin{subfigure}{0.325\textwidth}
    \begin{tikzpicture}
        \begin{axis}
            [width=\textwidth, height=\hei, xlabel=Iterations, ylabel=, tick label style={font=\normalsize},line width=1.5pt,legend style={draw=none, fill opacity=0.7, text opacity= 1,row sep=2pt, font=\normalsize, at={(0.95,0.5)}, anchor=east}, legend columns=1]
        \addplot+[mark=none,color=red] table [x=Time, y=phiAvg, col sep=comma] {pm-one-phase-back.csv};
        \addplot+[mark=none,color=blue] table [x=Time, y=phiAvg, col sep=comma] {pm-one-phase-forw.csv};
        \addplot+[mark=none,color=black] coordinates {(0,0.1) (682, 0.1)};
        \legend{Ref. conf., Forward, Given ($\overline\phi$)}
        \end{axis}
    \end{tikzpicture}
    \caption{Primal}
    \end{subfigure}
    \begin{subfigure}{0.325\textwidth}
    \begin{tikzpicture}
        \begin{axis}
            [width=\textwidth, height=\hei, xlabel=Iterations, ylabel=, tick label style={font=\normalsize},line width=1.5pt,legend style={draw=none, fill opacity=0.7, text opacity= 1,row sep=2pt, font=\normalsize, at={(0.95,0.5)}, anchor=east}, legend columns=1]
        \addplot+[mark=none,color=red] table [x=Time, y=phiAvg, col sep=comma] {pm-one-phase-mixed-p-back.csv};
        \addplot+[mark=none,color=blue] table [x=Time, y=phiAvg, col sep=comma] {pm-one-phase-mixed-p-forw.csv};
        \addplot+[mark=none,color=black] coordinates {(0,0.1) (682, 0.1)};
        \legend{Ref. conf., Forward, Given ($\overline\phi$)}
        \end{axis}
    \end{tikzpicture}
    \caption{Mixed $p$}
    \end{subfigure}
    \begin{subfigure}{0.325\textwidth}
    \begin{tikzpicture}
        \begin{axis}
            [width=\textwidth, height=\hei, xlabel=Iterations, ylabel=, tick label style={font=\normalsize},line width=1.5pt,legend style={draw=none, fill opacity=0.7, text opacity= 1,row sep=2pt, font=\normalsize, at={(0.95,0.5)}, anchor=east}, legend columns=1]
        \addplot+[mark=none,color=red] table [x=Time, y=phiAvg, col sep=comma] {pm-one-phase-mixed-u-back.csv};
        \addplot+[mark=none,color=blue] table [x=Time, y=phiAvg, col sep=comma] {pm-one-phase-mixed-u-forw.csv};
        \addplot+[mark=none,color=black] coordinates {(0,0.1) (841, 0.1)};
        \legend{Ref. conf., Forward, Given ($\overline\phi$)}
        \end{axis}
    \end{tikzpicture}
    \caption{Mixed $u$}
    \end{subfigure}

    \caption{Evolution of the average spatial porosity $\frac 1{|\Omega|}\int_\Omega\phi\,dx$ along the iterations of both the reference configuration and forward problems. We show the evolution obtained using the (a) primal formulation, (b) the mixed pressure formulation and (c) the mixed velocity formulation.}
    \label{fig:test-one-phase-mixed}
\end{figure}
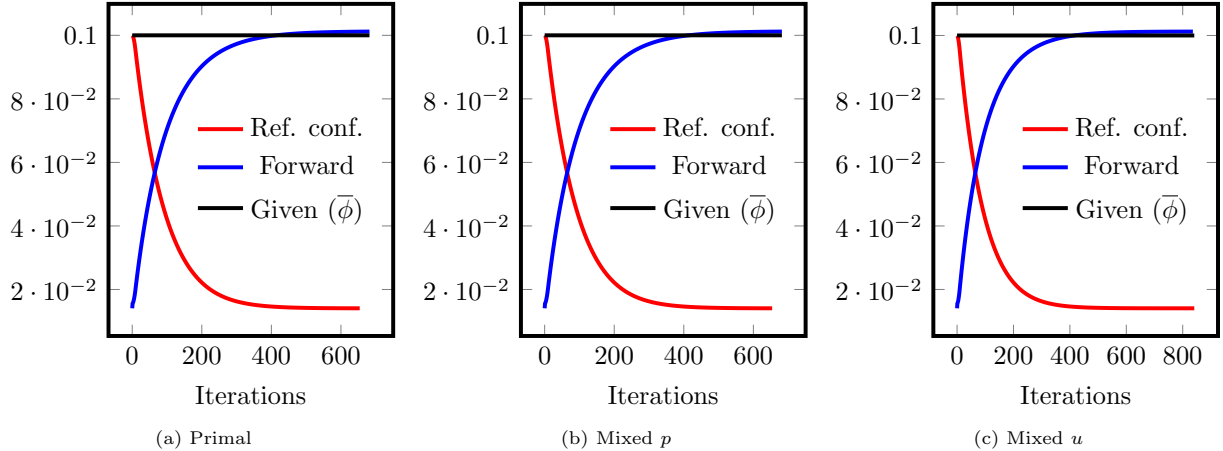

This test computes the reference configuration of the 2D square geometry as in Test~1, using the mixed models presented in Section~\ref{section:mixed}. The evolution of the average porosity in all three models is shown in Figure~\ref{fig:test-one-phase-mixed}, where all models yield the same solution. Indeed, the average reference porosity obtained by the primal, mixed in pressure, and mixed in velocity  models is given by $0.140767$, $0.140767$, and $0.140597$ respectively. We chose to prioritize stability in order to have an accurate comparison of the methods, so we used second order elements for the fluid velocity in the mixed velocity model. In addition, both primal and mixed pressure models took the same amount of iterations to converge (680), while the mixed velocity model took 842 iterations to converge.

\subsection{Test 5: Realistic perfusion model}
\begin{figure}
    \begin{subfigure}{0.3\textwidth}
        \includegraphics[width=\textwidth]{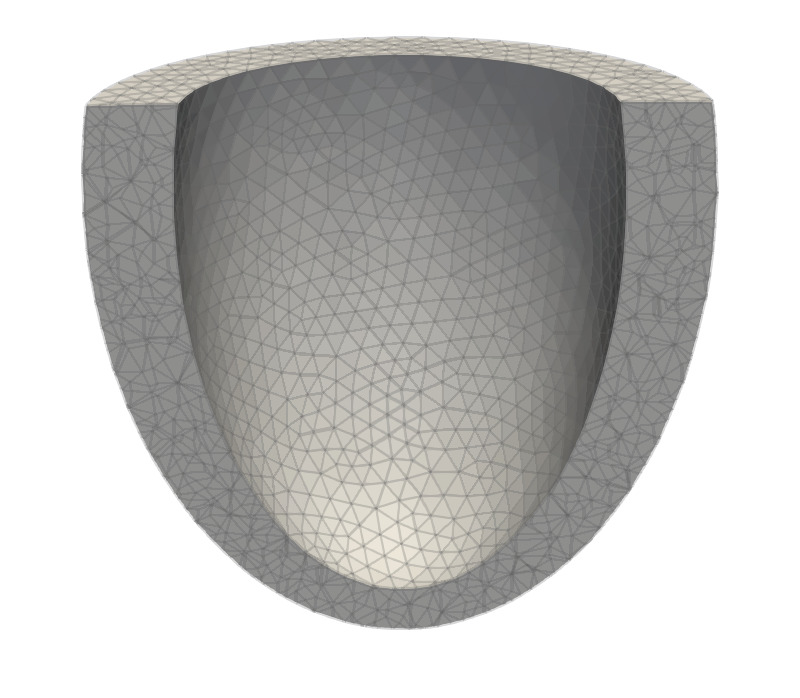}
        \caption{$\overline \phi=0.1$}
    \end{subfigure}
    \begin{subfigure}{0.3\textwidth}
        \includegraphics[width=\textwidth]{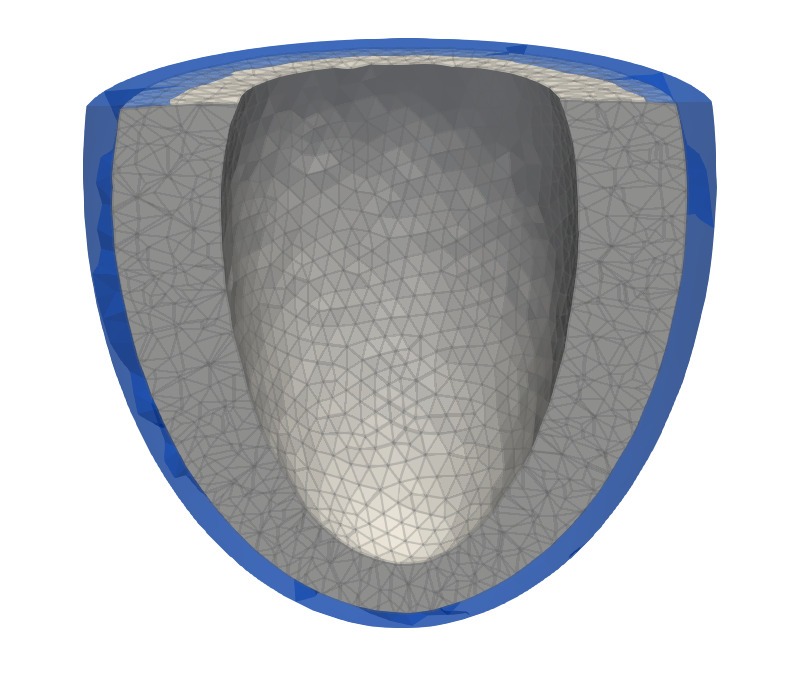}
        \caption{ $\phi_{h,0}=0.075$ }
    \end{subfigure}
    \begin{subfigure}{0.3\textwidth}
        \includegraphics[width=\textwidth]{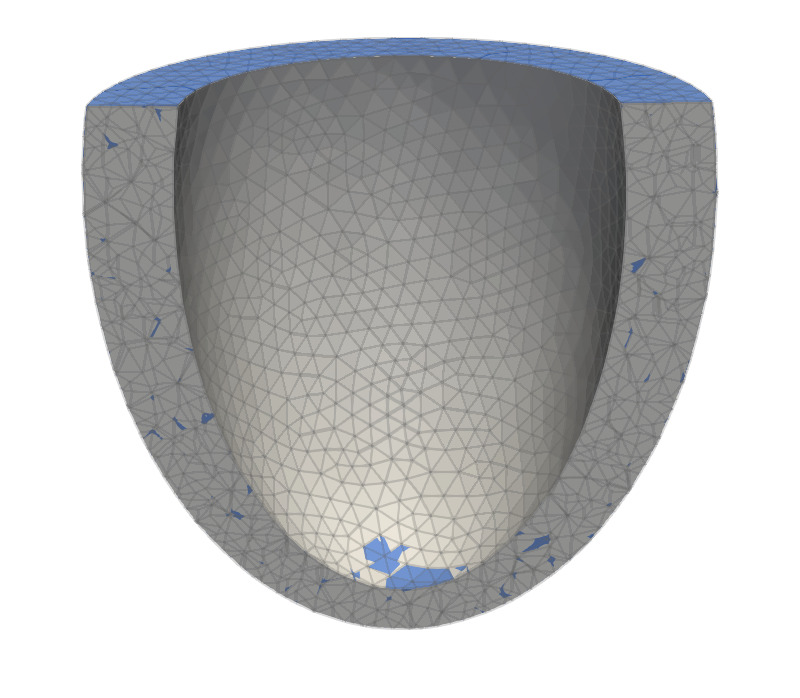}
        \caption{ $\phi_{h}=0.098$ }
    \end{subfigure}
    \caption{The average Eulerian steady steate porosities in the caption for the (a) given geometry, (b) computed reference configuration with the given geometry in a traslucent blue, and (c) the deformed geometry from the given loads with the given geometry in a traslucent blue. The 'h' subscript shows that the function is a FEM approximation.}
    \label{fig:lv-solution}
\end{figure}

In this test, we compute the reference configuration associated to a perfusion model of an idealized left ventricular geometry. For this, the modifications to the model considered are: (i) boundary conditions for the mechanics are divided into an endocardial pressure of $1.5\,\texttt{kPa}$ (physiological end-diastolic pressure) and friction due to the pericardium on the epicardium as in \cite{pfallerimportance}, (ii) the source term $\theta$ considers both an arterial source and a venous sink of $5\,\texttt{kPa}$ and $1\,\texttt{kPa}$ respectively, (iii) the passive material is the anisotropic one from \cite{Usyk2000}, with the fiber, sheet and normal directions computed using an accurate $H^1$ conforming approximation proposed in \cite{barnafi2023modeling}, and (iv) we replace the solid incompressibility with a quasi-incompressible law of the same magnitude of the geometry quasi-incompressibility as in \cite{barnafi2022perfusion}. The geometry, the computed reference configuration and the forward solution are displayed in Figure~\ref{fig:lv-solution}, together with their corresponding reference porosities in the caption. The solution displays a very good match with respect to the given geometry, and the average stationary porosity yields a small error of a 2\% (0.098 instead of 0.1) due to numerical accuracy. Additionally, the average stationary porosity obtained with the primal model is 0.094, showing that not using a consistent mixed formulation can yield in this simple test an additional 4\% error. We expect this number to increase in more complex scenarios.

\begin{figure}
    \newcommand{\hei}{5cm}
    \centering
    \begin{subfigure}{0.49\textwidth}
    \begin{tikzpicture}
        \begin{axis}
            [width=\textwidth, height=\hei, xlabel=Time, ylabel=, tick label style={font=\normalsize},line width=1.3pt,legend style={draw=none, fill opacity=0.7, text opacity= 1,row sep=2pt, font=\normalsize, at={(0.95,0.5)}, anchor=east}, legend columns=2]
        \addplot+[mark=none,color=red] coordinates {(0.1,162) (0.2,169) (0.3,171) (0.4,173) (0.5,175) (0.6,177) (0.7,179) (0.8,181) (0.9,184) (1.0,187)};
        \addplot+[mark=o,color=green!70!black] coordinates {(0.1,19) (0.2,21) (0.3,23) (0.4,25) (0.5,27) (0.6,29) (0.7,31) (0.8,33) (0.9,37) (1.0,39)};
        \addplot+[mark=triangle,color=blue!70!white] coordinates {(0.1,17) (0.2,19) (0.3,21) (0.4,23) (0.5,25) (0.6,27) (0.7,29) (0.8,31) (0.9,34) (1.0,37)};
        \addplot+[mark=square,color=blue!30!white] coordinates {(0.1,36) (0.2,38) (0.3,40) (0.4,42) (0.5,44) (0.6,46) (0.7,48) (0.8,50) (0.9,53) (1.0,56)};
        \legend{$AA(0)$, $AA(1)$, $AA(2)$, $AA(5)$}
        \end{axis}
    \end{tikzpicture}
    \caption{Ref. conf.}
    \end{subfigure}
    \begin{subfigure}{0.49\textwidth}
    \begin{tikzpicture}
        \begin{axis}
            [width=\textwidth, height=\hei, xlabel=Time, ylabel=, tick label style={font=\normalsize},line width=1.3pt,legend style={draw=none, fill opacity=0.7, text opacity= 1,row sep=2pt, font=\normalsize, at={(0.95,0.5)}, anchor=east}, legend columns=2]
        \addplot+[mark=none,color=red] coordinates {(0.1,38) (0.2,43) (0.3,45) (0.4,46) (0.5,47) (0.6,48) (0.7,49) (0.8,50) (0.9,51) (1.0,52)};
        \addplot+[mark=o,color=green!70!black] coordinates {(0.1,13) (0.2,14) (0.3,15) (0.4,16) (0.5,17) (0.6,18) (0.7,19) (0.8,20) (0.9,21) (1.0,22)};
        \addplot+[mark=triangle,color=blue!70!white] coordinates {(0.1,10) (0.2,11) (0.3,12) (0.4,13) (0.5,14) (0.6,15) (0.7,16) (0.8,17) (0.9,18) (1.0,19)};
        \addplot+[mark=square,color=blue!30!white] coordinates {(0.1,10) (0.2,11) (0.3,12) (0.4,13) (0.5,14) (0.6,15) (0.7,16) (0.8,17) (0.9,18) (1.0,19)};
        \legend{$AA(0)$, $AA(1)$, $AA(2)$, $AA(5)$}
        \end{axis}
    \end{tikzpicture}
    \caption{Forward}
    \end{subfigure}

    \caption{Number of accumulated time iterations incurred to solve the (a) reference configuration and (b) forward problems at each ramp step. In (b), the $AA(2)$ and $AA(5)$ lines overlap.}
    \label{fig:lv-iters}
\end{figure}
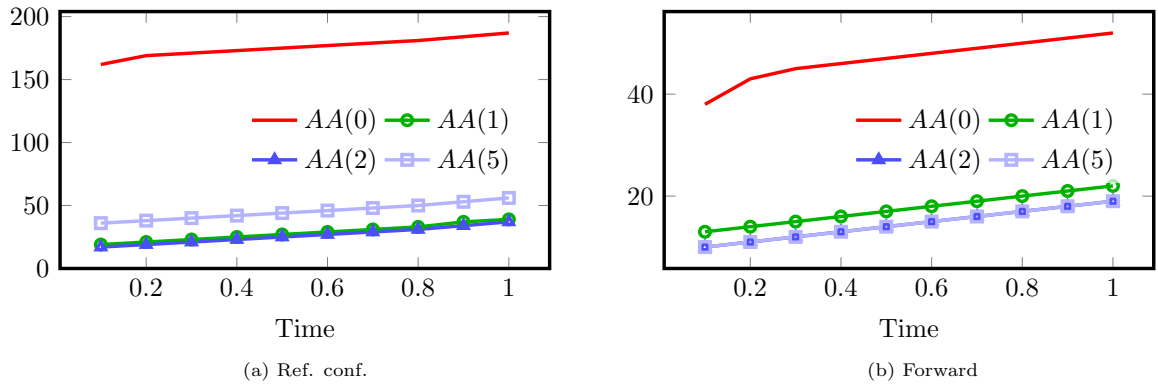

We note that this test was significantly more difficult than the others, as it would frequently result in diverging nonlinear iterations. To circumvent this difficulty, instead of increasing our ramp slowly and then computing a stationary state, we computed the stationary state for each ramp level before increasing it. A natural improvement for this was Anderson acceleration, so in Figure~\ref{fig:lv-iters} we show the number of fixed point iterations accumulated at each ramp level for 1, 2, and 5 previous vectors, and in Table~\ref{table:lv-iters} we show the total number of iterations (rightmost value at Figure~\ref{fig:lv-iters}). We note that iterations are reduced by up to roughly a 63\% and an 80\% of the ones obtained without acceleration ($AA(0)$) for the forward and reference configuration problems respectively, and interestingly, the optimal number of previous vectors depends on the formulation. Still, considering also the previous tests, we believe that a depth of 2 is a good general choice.

\begin{table}
    \centering
    \begin{tabular}{r | r r r r}
    \toprule & $AA(0)$ & $AA(1)$ & $AA(2)$ & $AA(5)$ \\\midrule
    Ref. conf. & 187     & 39      & 37      & 56 \\\midrule
    Forward  & 52      & 22      & 19      & 19 \\\bottomrule
    \end{tabular}
    \caption{Total number of fixed point iterations incurred when solving both reference configuration and forward models in the left ventricle test.}
    \label{table:lv-iters}
\end{table}

\section{Discussion}\label{section:discussion}
This work proposes a novel model for computing the reference configuration of a fully nonlinear poroelastic material, where the main ingredient is the identification of the given porous solution as the steady state configuration of a porous media problem with unknown initial conditions. The proposed model is consistent with fundamental thermodynamics principles, which renders it applicable for all applications where nonlinear poroelasticity is relevant. We have tested the model using parameters and laws coming from cardiac applications, which involve the combination of logarithmic and exponential nonlinearities that render the model computationally challenging.

One interesting discovery in this work is the \emph{primal inconsistency} of the primal formulation, which can be avoided using mixed formulations. We proposed two approaches, one that yields a natural energy norm and another one that yields a saddle point problem. The saddle point formulation is much more challenging, but it can be useful when considering more general fluid models. Our main interest is to consider Darcy-Brinkman models for cardiac ablation, which has been shown to be relevant in applications \cite{wongchadakul2023tissue}. However solved, the computation of a stationary state is an expensive problem, which results in simulations with hundreds of iterations. To circumvent this difficulty, we have shown that Anderson acceleration yields an extremely successful strategy to reduce iteration counts, with reductions of up to an 80\%.

Our future work mainly focuses on two topics: the formulation of robust reference configuration models for multi-phase porous media, which we observed to be much more difficult to solve, and the use of the proposed models to accurately model the interplay between blood flow and deformation during cardiac ablation as in \cite{petras2019computational}, by extending the large deformations model proposed in \cite{molinari2022transversely}.

\section{Acknowledgments}
NB has been funded by the ANID Postdoctoral 3230326, by CMM BASAL FB2100005 and partially supported by Johannes Kepler University.  AP and LGG were partially supported by the State of Upper Austria.

\bibliography{main} 
\bibliographystyle{alpha}
\end{document}